\documentclass[a4paper,leqno]{amsart}

\usepackage{verbatim}

\usepackage{amsmath,enumerate,rotate}
\usepackage{amssymb,latexsym}
\usepackage{epsfig}
\usepackage{graphicx}
\usepackage[english]{babel}

\newtheorem{definition}{Definition}[section]
\newtheorem{theorem}[definition]{Theorem}
\newtheorem{lemma}[definition]{Lemma}
\newtheorem{proposition}[definition]{Proposition}

\newtheorem{remark}{Remark}

\def\X{X}

\def \fine {\diamondsuit}

\def \lf {\left}
\def \rg {\right}

\def\a {\alpha}
\def\t {\tau}
\def\th {\theta}

\def\s {\sigma}

\def\eps {\epsilon}

\def\b {\beta}

\def\lm {\lambda}

\def\J {\mathbb{I}}

\def\B {\mathbb{B}}

\def\N {\mathbb{N}}

\def\RE {\mathbb{R}}

\def\N {\mathbb{N}}

\def\CG {\mathcal{G}}

\def\CS {\mathcal{S}}
\def\CE {\mathcal{E}}

\begin{document}

\numberwithin{equation}{section}
\def \X {\mathcal{X}}
\def \CX {\mathcal{X}}
\def\B {\mathcal{B}}

\title[Metropolis algorithm and equienergy sampling]{Metropolis algorithm and  equienergy sampling for two mean field spin systems}


\author{Federico Bassetti and Fabrizio Leisen }

\address{Universit\`a degli Studi di  Pavia, Dipartimento di Matematica, via Ferrata 1, 27100 Pavia, Italy}

\address{Universit\`a dell'Insubria, Dipartimento di
  economia, via monte generoso 71, 21100 Varese , Italy}

\email{federico.bassetti@unipv.it}
\email{leisen.fabrizio@unimore.it}

\begin{abstract}
{In this  paper we study   the Metropolis  algorithm
in connection with two mean--field spin systems,
 the so called mean--field Ising model and the
Blume--Emery--Griffiths model. In both this examples
 the naive choice of proposal chain gives
rise, for some parameters, to a
slowly mixing Metropolis chain,  that is a chain whose
 spectral gap decreases
exponentially fast (in the
dimension $N$ of the problem).
Here we show how a slight variant in the proposal chain can
avoid this problem,
keeping the  mean computational cost similar to the cost of the
usual Metropolis. More precisely we prove that,
with  a suitable variant in the proposal,
 the Metropolis chain has a spectral  gap which decreases  polynomially
 in $1/N$.
Using some symmetry structure of the energy, the method rests on allowing
appropriate jumps within the energy level of the starting state, and
it is strictly connected to both
the {\it small world Markov chains} of \cite{guan1,guan2}
and to the {\it equi-energy sampling} of \cite{Kou2006} and
\cite{MadrasPiccioni}.
}

\end{abstract}


\maketitle

\section{Introduction.}



The Metropolis algorithm,  introduced in
\cite{Metropolis} and later generalized in \cite{Hastings},
 is currently
(together with  other  Monte Carlo Markov Chain methods)
one of the most used
simulation techniques both in statistics and in physics.
See, among others,
\cite{Peskun2,Peskun,Sokal,pippo,Rubinstein,Robert,Liu2,DSC2}.

In a finite setting the Metropolis algorithm can be described as
follows.
Suppose that, given a probability $\pi(x)$ on a finite set $\CX$,
want to approximate
\begin{equation}\label{1.0}
\mu=\sum_x f(x)\pi(x),
\end{equation}
for $f:\CX \to \RE$.
As a first step,
take a reversible Markov chain $K(x,y)$ (the
proposal chain) on
$\CX$ 
and
change its output in order  to have a new chain with stationary
distribution $\pi$.
This can be achieved  by
constructing a new ($\pi$--reversible) chain
\begin{equation}\label{MChain}
M(x,y) = \left \{
\begin{array}{lll}
K(x,y)A(x,y) &\quad  x \, \not = y
 \\
K(x,x)+\sum_{z \not = x }
K(x,z)(1-A(x,z)) & \quad x=y \,\\
\end{array}
\right .
\end{equation}
where $A(x,y):=\min(\frac{\pi(y)K(y,x)}{\pi(x)K(x,y)},1) $.
Then,
the {\it metropolis estimate} of $\mu$ is given by
\begin{equation}\label{1.1}
\hat \mu_n= \frac{1}{n} \sum_{i=1}^n f(Y_i),
\end{equation}
where $Y_0$ is generated from some initial distribution
$\pi_0$ and  $Y_1,\dots,Y_n$ from $M(x,y)$.

It is clear that, from a computational  point of view,
  the speed of convergence
to the stationary distribution and the (asymptotic) variance of the estimate
are two very important features of the Markov chain $M$.

It is well-known that in some situation a Markov chain can converge
very slowly to its stationary distribution and, moreover, that  the
asymptotic variance of the estimate (\ref{1.1})  can be much bigger
than the variance of $f$, i.e. $Var_\pi(f):=\sum_x(f(x)-\mu)^2\pi(x)$, which is
equal to the asymptotic variance of the crude Montecarlo estimator.
In these cases (\ref{1.1}) turns out to be a very
inefficient estimate of $\mu$.


For the Metropolis chain
a classical  situation in which the
convergence is slow (and the  variance big) is when
the target distribution $\pi$ has  many peaks and $K$
is somehow too ``local''.


This is well known in statistical physics, where, typically,
a distribution of a system with
energy function $h$ and in thermal
equilibrium at temperature $T$ is described by the Gibbs distribution
\[
\pi_{h,T}(x)=\exp\{ -h(x)/T \}Z_T^{-1}
\]
with $Z_T=\sum_x \exp\{ -h(x)/T \}$. In point of fact,
the Metropolis algorithm has
been proposed in \cite{Metropolis} to compute average
with respect to such distributions. Indeed, if $h$ is nice, the Metropolis
algorithm
is very
efficient, but it can perform very poorly if the energy has many
local minima separated by high barriers that cannot be crossed by
the proposal moves $K$. This problem can be  bypassed, for specific energy,
designing appropriate moves that have higher chance to cut across the
energy barrier (see, e.g, \cite{Caracciolo1990,Caracciolo1991}),
or constructing clever alternative approaches to the problem,
for instance using a reparametrization of the problem
(see, e.g.,  \cite{Gelfand1995,Gelman1996}) or using
auxiliary variables (see, e.g., \cite{Swendsen87,Sokal88,Besag1993,Mira2002}).
A different kind of
solution  has been proposed in  \cite{Geyer}
and in \cite{MarinariParisi92}
 by introducing  the so
called {\it simulated tempering}, which essentially means that $T$ is changed
(stochastically or not) to flatten $h$. A remarkable variant of these
methods is the {\it parallel tempering}, see, for instance, \cite{nemoto}.
More recently new algorithms
based on the so called
{\it equi--energy levels sampling} have been proposed (see
\cite{MadrasPiccioni} and  \cite{Kou2006}). In particular, the
  algorithm proposed
in \cite{Kou2006}
relies on  the so--called equi-energy
jump, which enables the chain to reach regions of the sample space
with energy close to the one of the starting state,
 but that  may be separated
by steep energy barriers.  In point of fact, even if, according to
some simulations,
 the method seems to be efficient nothing has been formally proved.
Finally, let us mention a recent algorithm,
called {\it small world Markov chains} (see \cite{guan1,guan2}),
 that combine a local chain with
long jumps.
In these papers, it has been shown that  a simple modification
of the proposal mechanism results in faster convergence of the chain.
That mechanism,
which is based on an idea from the field of "small-world"
networks, amounts to adding occasional "wild" proposals to any
local proposal scheme.

In the present paper we study two simple
examples: the so called mean field Ising model and the mean field
Blume--Emery--Griffiths model. As for the former, it is well-known
that  the usual choice of $K$ gives rise, for low temperature, to a
slowly mixing Metropolis chain  (see, e.g., \cite{MadrasPiccioni}).
Here we show that a slight variant in the proposal chain can completely
solve this problem, keeping the  mean computational cost similar to
the cost of the usual Metropolis. The idea again rests on allowing
appropriate jumps in the same energy level of the starting state. As for
the Blume--Emery--Griffyths mean--field
model, we first show that there is a critical region of the parameters space for which
the naive Metropolis chain is slowly mixing. Then we show how one can modify the proposal
chain in order to obtain a better mixing for the Metropolis chain.
The present paper should be intended as a further
 step in the direction of a better mathematical
understanding of both { small world Markov chains}
and {equi-energy sampling}.

The rest of the paper is organized as follows.
In Section \ref{S:generalization} some general considerations are given.
In Section \ref{S:start} some basic tools concerning Markov chain,
which will
be used in the paper,
are reviewed.
Section \ref{S:warmingup} contains a warming up example.
In Section \ref{S:meanfield} the mean field Ising model is treated, while
 Section \ref{S:BEG} deals with the more complex case of the mean field
Blume-Emery-Griffiths model.
 All the proofs are deferred to the Appendix.

\section{A general strategy}\label{S:generalization}
In an abstract setting, what we
shall do  in the next examples can be summarized as follows.
Let $\CG$ be a group acting on $\CX$ for which
\begin{equation}\label{symmetry}
\pi(x)=\pi(g(x)) \qquad \forall \,\, x\,\,\in \CX, \,\, \forall \,\, g\,\,\in \CG.
\end{equation}
For every $x$ in $\X$
let \(
O_x:=\{y=g(x): g \in \CG \}
\) be the orbit of $x$ (of course if $y$ belongs to $O_x$ then $O_x=O_y$).

Assume now that we have a reversible Markov chain $K_E(x,y)$ (the
proposal) on $\CX$ and suppose that the Metropolis chain  $M_E$
with proposal $K_E$ is slowly mixing (see next section for more details).
To speed up the mixing one can try to exploit (\ref{symmetry}) by taking a proposal
of the following form:
\begin{equation}\label{smallw}
K_\eps(x,y)=\eps K_E(x,y)+(1-\eps)K_\CG(x,y)
\end{equation}
where
\[
K_\CG(x,y)= \sum_{z \in O_x }q_x(z)
\J_{z}(y),
\]
$0<q_x(z)<1$ and $\sum_{z \in O_x }q_x(z) =1$.

In point of fact, usually $K_E$ is ``local''; for instance frequently
\[
K_E(x,y)=0
\]
whenever $y \not=x$ belongs to $O_x$, hence with $K_\CG$ we are adding
``long'' jumps to the chain.
Moreover, note that if $K_E$ is such that $K_E(x,g(x))=K_E(g(x),x)$,
for every $x$ in $\CX$ and $g$ in $\CG$, then the Metropolis always accepts the move
$x \to g(x)$ and
\[
M(x,g(x))=\eps K_E(x,g(x))+ (1-\eps) q_x(g(x)).
\]
In particular this holds when $K_E$ is symmetric.

The heuristics under (\ref{smallw}) is to combining {small world Markov chains}
and {equi-energy sampling}.

Before presenting some examples in which one can actually
improve the performances of the Metropolis chain using this idea,
we collect in the next section
some useful facts concerning Markov chains.

\section{Preliminaries}\label{S:start}
Let $P(x,y)$ be a reversible and ergodic
 Markov chain on the finite set $\CX$ with (unique) stationary distribution
$p(x)$. Thus, $p(x)P(x,y)=p(y)P(y,x)$.
Let $L^2(p)=\{f:\CX \to \RE  \}$ with
$<f,g>_p=E_p(fg)=\sum_x f(x) g(x) p(x)$. Reversibility
is equivalent to $P:L^2 \to L^2$ being self--adjoint.
Here $Pf(x)=\sum_y f(y) P(x,y)$.
The spectral theorem implies that $P$ has real eigenvalues
$1=\lambda_0(P) > \lambda_1(P) \geq \lambda_2(P) \geq \dots \geq \lambda_{|\CX|-1}(P) > -1$
with orthonormal basis of eigen--functions
$\psi_i:\CX \to \RE$ ($P \psi_i(x)=\lambda_i \psi_i(x)$,
$<\psi_i,\psi_j>_p=\delta_{ij}$).

\subsection{Spectral gap, variance and speed of convergence}
A very important quantity related to the eigenvalues is the spectral gap, defined by
\[
Gap (P) =1 -\max\{\lambda_1,|\lambda_{|\X|-1}| \}.
\]
It turns out that the spectral gap is a good index to measure the mixing of a chain.
To better understand this point, assume that $f$ belongs  to $L^2(p)$
and write
$f(x)=\sum_{i \geq 0}a_i \psi_i(x)$ {\rm(}with
$a_i=<f,\psi_i>_p${\rm)}.
Now let $Y_0$ be chosen form some distribution $p_0$
and $Y_1, \dots, Y_n$ be a realization of the $P(x,y)$ chain, then
\[
\hat \mu_n= \frac{1}{n} \sum_{i=1}^n f(Y_i)
\]
has asymptotic variance
given by
\[
AVar(f,p,P):=\lim_{n \to +\infty} n \cdot Var(\hat \mu_n)=
\sum_{k \geq 1} |a_k|^2 \frac{1+\lambda_k}{1-\lambda_k}.
\]
See, for instance, Theorem 6.5 in Chapter 6 of \cite{Bremaud}.
From the last expression,
the classical inequality
\begin{equation}\label{avarineq}
AVar(f,p,P) \leq \frac{2}{1-\lambda_1} Var_{p}(f),
\end{equation}
  follows easily.
The last inequality is the
usual way of relating spectral gap to asymptotic variance and, hence,
to the efficency of a chain.

The spectral gap is very important also
to give bounds on the speed of convergence to the stationary distribution.
For example, if
$\|\cdot\|_{TV}$ denotes the total variation norm,
one has
\[
 \|\delta_x P^k-p\|_{TV}^2 =\left(\sup_{A \subset \X}|P^k(x,A)-p(x)|\right)^2
\leq \frac{1-p(x)}{4p(x)}(\max\{\lambda_1,|\lambda_{|\X|-1}|\})^{2k}
\]
See, e.g., Proposition 3 in \cite{DiacStrook}.
Another classical bound is
\[
\|p_0P^k/p-1\|_{2,p} \leq Gap(P^k) \|p_0/p -1\|_{2,p}
\]
valid for every probability $p_0$. See, for instance, \cite{Sokal}.

Roughly speaking one can say that a sequence of Markov chains defined on a sequence of state space
$\CX_N$ is slowly mixing (in the dimension of the problem $N$)
if the spectral gap decreases exponentially fast in  $N$.

\subsection{Cheeger's inequality}
As already recalled, problems of slowly mixing
typically occur when
 $\pi$ has two or more peaks and the chain
$K$ can
only move in a neighborhood of
 the starting peak.
Usually this phenomenon is called bottleneck.
A powerful tool to detect the presence of a bottleneck
 is the conductance and the related Cheeger's inequality.
Recall that the conductance of a chain $P$  with stationary
distribution
$p$ is defined by
\[
h=h(p,P):=\inf_{A \,\, : p(A) \leq  \frac{1}{2}  }  \,\,\, \frac{1}{p(A)}
\sum_{x \in A,  y \in A^c} p(x) P(x,y),
\]
and the well-known Cheeger's inequality is
\begin{equation}\label{cheger}
1-2 h \leq \lambda_1(P) \leq 1- \frac{h^2}{2}.
\end{equation}
See, for instance, \cite{Bremaud,sinclair1993,DiacStrook}.
Note that, since $P$ is reversible,
\begin{equation}\label{leqcond}
h \leq \frac{1}{p(A)} \sum_{x \in A}\sum_{y \in A^c} p(x) P(x,y)= \frac{1}{p(A)}
\sum_{x \in A}\sum_{y \in A^c} p(y) P(y,x)
\end{equation}
for every $A$ such that $p(A) \leq 1/2$.

\subsection{Chain decomposition theorem}\label{chaindec}
In this subsection we briefly describe a useful
technique to obtain bounds on the spectral gap: the so called chain
decomposition technique. Following \cite{guan2}
assume that $A_1,\dots,A_m$ is a partition
of $\X$. Moreover, for each $i=1,\dots,m$, define a new Markov chain on $A_i$
by setting
\[
P_{A_i}(x,y):=P(x,y)+\J_{x}(y)\left (\sum_{z  \in A_i^c }P(x,z) \right) \qquad
(x,y \in A_i).
\]
$P_{A_i}$ is a reversible chain  on the state space $A_i$ with respect to
the probability measure
\[
p_i(x):= p(x)/p(A_i).
\]
The movement of the original chain among the ``pieces''
$A_1,\dots,A_m$ can be
described
by a Markov chain with state space $\{1,\dots,m\}$ and transition
probabilities
\[
P_H(i,j):=\frac{1}{2 p(A_i)} \sum_{x \in A_i, y \in A_j} P(x,y)p(x)
\]
for $i \not = j$ and
\[
P_H(i,i):=1 - \sum_{j \not = i} P_H(i,j),
\]
which is reversible with stationary distribution
\[
\bar p(i):=p(A_i).
\]
A variant of a result of Caracciolo, Pelisetto and Sokal (published in \cite{MadrasRandall}),
 states that
\begin{equation}\label{decLB}
Gap(P) \geq \frac{1}{2} Gap(P_H)
\left ( \min_{i=1,\dots,m} Gap(P_{A_i}) \right)
\end{equation}
holds true, see Theorem 2.2 in  \cite{guan2}.
Other results about
chain decompositions can be found, for instance, in  \cite{jerrum2004}.

In the next very simple example we shall show how this
technique can be used,
starting from a slowly mixing chain, to suggest how to modify
the proposal chain in order to obtain a fast mixing chain.

\section{Warming up example}\label{S:warmingup}
Set $\X=\{-N,-N+1,\dots,0,1,\dots, N \}$
and define a probability measure on $\X$ by
\[
\pi(x)= \frac{(\th-1)\th^{|x|}}{2\th^{N+1}+1-\th},
\]
$\th$ being a given parameter bigger than $1$.
Here we can consider $\CG=\{+1,-1\}$ (with group operation given by
the usual product) acting on $\CX$ by $g(x)=g x$, hence
$O_x=\{x,-x\}$.

Now let $K_E$ be a chain defined by
\[
\begin{split}
&K_E(x,x+1) =1/2  \qquad x \not = N \\
&K_E(x,x-1) =1/2  \qquad x \not = -N \\
&K_E(N,N)=K_E(-N,-N)=1/2 \\
&K_E(x,y)=0 \quad \text{otherwise}
\end{split}
\]
and denote by $M_E$ the Metropolis chain with stationary
distribution $\pi$  derived by $K_E$. It is clear that in this case
$K_E(x,y)=0$ whenever $y$ belongs to $O_x$. In this example
it is very easy to bound the conductance on $M_E$, indeed, taking
$A=\{-N,\dots,-1\}$, by (\ref{leqcond}), it follows that
\[
h(\pi,M_E) \leq \frac{\pi(0)}{1-\pi(0)}.
\]
Hence,
\[
h(\pi,M_E) \leq C \th^{-N},
\]
and then (\ref{cheger}) yields
\[
1-\lambda_1 \leq 2C \th^{-N}.
\]
This means that, if $f$ is such that $a_1 \not =0$ and $\th>1$, then
the asymptotic variance of $f$ blows up exponentially fast, indeed
\[
AVar(f,\pi,M_E) \geq 2C e^{\log(\th) N}.
\]

Now, instead
  of $K_E$ consider
\[
K_\eps(x,y)=(1-\eps)K_E(x,y)+ \eps\J_{\{-x\}}(y)
\]
and let $M^{(\eps)}$ be the Metropolis chain derived by $K_\eps$.
Decompose $\X$ as follows
\[
\X=A_1 \cup A_2 \dots \cup A_N
\]
with $A_1=\{-1,0,1\}$ and $A_i=\{x \in \X:
|x|=i\}$, for $i>1$. Moreover let
\[
\bar \pi(i)=\pi(A_i)=
\left \{
\begin{array}{ll}
(2\th+1)/Z & \text{for $i=1$}\\
2 \th^i/Z  &  \text{for $i>1$}
\end{array}
\right .
\]
where
\[
Z=\frac{2\th^{N+1}+1-\th}{(\th-1)}
\]
and set
\[
M^{(\eps)}_H(i,j)=\frac{1}{2\pi(A_i)}\sum_{l \in A_i,m \in
  A_j}M^{(\eps)}(l,m)\pi(l),
\qquad M^{(\eps)}_H(i,i)=1-\sum_{j\not=i}M^{(\eps)}_H(i,j).
\]
For $i\not=1,N$, one has
\[
M^{(\eps)}_H(i,i+1)= \frac{1}{2\pi(A_i)}[M^{(\eps)}(i,i+1)\pi(i)+
M^{(\eps)}(-i,-i-1)\pi(-i)]
\]
and, since $\pi(i)=\pi(-i)$ and $\pi(i+1) \geq \pi(i)$
\[
M^{(\eps)}_H(i,i+1)=\frac{1-\eps}{4}.
\]
In the same way it is easy to see that
\[
\begin{split}
&M^{(\eps)}_H(i,i-1)=\frac{1-\eps}{4\th}, \qquad i\not=1,N\\
&M^{(\eps)}_H(i,i)=1-\frac{1-\eps}{4}(1+\th^{-1})  \qquad i\not=1,N\\
&M^{(\eps)}_H(N,N-1)=\frac{1-\eps}{4\th} \qquad M^{(\eps)}_H(N,N)=1-\frac{1-\eps}{4\th} \\
&M^{(\eps)}_H(1,2)=\frac{1-\eps}{4(1+1/(2\th))} \qquad
M^{(\eps)}_H(1,1)=1-\frac{1-\eps}{4(1+1/(2\th))}. \\
\end{split}
\]
Moreover, for every $i\not =1$, $M^{(\eps)}_{A_i}$ in matrix form is given by
\begin{equation*}
\left (
\begin{array}{ll}
1-\eps & \eps \\
\eps &  1-\eps
\end{array}
\right ),
\end{equation*}
and hence
\[
Gap(M^{(\eps)}_{A_i})=1-|1-2\eps|.
\]
While $M^{(\eps)}_{A_1}$
is given by
\begin{equation*}
\left (
\begin{array}{ccc}
(2\th-1)(1-\eps)/(2\th)  & (1-\eps)/(2\th) & \eps  \\
(1-\eps)/2 & \eps &  (1-\eps)/2 \\
\eps &  (1-\eps)/(2\th) & (2\th-1)(1-\eps)/(2\th)
\end{array}
\right )
\end{equation*}
and hence
\[
Gap(M^{(\eps)}_{A_i})=k(\th,\eps)>-1.
\]
Moreover, since
\[
\begin{split}
\min & \lf [\min_{i\not =1,N} (M^{(\eps)}_H(i,i\pm 1)),
M^{(\eps)}_H(1,2),M^{(\eps)}_H(N,N-1) \rg ] \\
&
 \geq \min \lf [(1-\eps)/(4\th),\frac{1-\eps}{4(1+1/(2\th))} \rg ]=:m(\eps,\th) >0 \\
 \end{split}
\]
and $\bar \pi(i) \leq 3 \bar \pi(j)$ for every $i<j$,
Lemma \ref{lemmaBD} in the appendix yields that
\[
1-\lambda_1(M^{(\eps)}_H) \geq \frac{m(\eps,\th)}{3N^2}.
\]
In the same way, since  $M^{(\eps)}_H(i,i+ 1) +
M^{(\eps)}_H(i,i- 1)
\leq (1-\eps)M(\th)/4$, with $M(\th)=\max(1+\th^{-1},2\th/(2\th+1))
\leq 2$,
inequality
(\ref{lowerboundBDchain}) in the Appendix yields that
\[
\lambda_{N-1}(M^{(\eps)}_H) \geq 1 -\frac{1-\eps}{2} \geq \frac{1+\eps}{2}.
\]
Hence
\[
Gap(M^{(\eps)}_H) \geq \frac{m(\eps,\th)}{3N^2}
\]
and (\ref{decLB}) yield
\[
Gap(M^{(\eps)}) \geq \frac{h(\th,\eps)}{N^2}
\]
for a suitable $h$.
This shows that $M^{(\eps)}$ is fast mixing for every $\eps>0$ and
for every $\th>1$ while
$M_E$ is slowly mixing for every $\th>1$.

\section{The mean field Ising model}\label{S:meanfield}
Let $\X=\{-1,1 \}^N$, $N$ being an even integer.
For every $\beta > 0$ let
$\pi=\pi_{\beta,N}$ be a probability on  $\X$ defined by
\[
\pi(x)=\pi_{\beta,N}(x):= \exp \left \{ \beta
\frac{S_N^2(x)}{2N} \right \}Z^{-1}_N(\beta)  \qquad (x \in \X)
\]
where
\[
Z_N(\beta)=Z_N:=\sum_{x \in \X} \exp \left \{ \beta
\frac{S_N^2(x)}{2N} \right \}
\]
is the normalization constant (``partition function'') and
\[
S_N(x):= \sum_{i=1}^N x_i \qquad x=(x_1,\dots,x_N).
\]
This is the so called {\it mean field Ising model}, or {\it Curie-Weiss model}, in which every particle
$i$, with spin $x_i$, interacts equally with every other particle.
It is probably the most simple but also the most studied example of spin system on a complete graph.
The usual Metropolis algorithm uses  as proposal chain
\[
K_E(x,y)=\frac{1}{N} \sum_{j=1}^N \J_{\{x^{(j)} \}}(y)
\]
where $x^{(j)}$ denotes the vector $(x_1,\dots,-x_j,\dots,x_N)$.
It has been proved in \cite{MadrasPiccioni}
that, whenever $\b >1$,
\[
1-\lambda_1 \leq C e^{-D^2 N}
\]
where $\lambda_1$ is the first eigenvalues
smaller than $1$ of  the Metropolis chain $M_E$ derived $K_E$.
This yields that
the variance of an estimator obtained from this Metropolis algorithm
can blow up exponentially fast in $N$.

The aim of this section is to show how one can construct a different
Metropolis
chain avoiding this problem. In the notation of Section \ref{S:generalization},
we consider
\[
\CG= \CS_N \times \{+1, -1 \}
\]
($\CS_N$ being the symmetric group of order $N$)
and we define the action of $\CG$ on $\X=\{-1,1\}^N$ by
\[
g(x)=(e\cdot x_{\s(1)},\dots,e\cdot x_{\s(N)}) \qquad g=(\s,e).
\]
In order to introduce a new proposal,  it is useful to
write $\X$ as the union of its ``energy sets'',
 that
is
\[
\X=\X_0 \cup \X_2 \cup \X_4 \cup \dots \cup \X_N
\]
where
\[
\X_i:=\{x \in \X : |S_N(x)|=i \} \qquad (i=0,2,\dots,N).
\]
Note that energy takes only even values and that
$O_x=\X_{|S_N(x)|}$.
Moreover, for $i \not = 0$, set
\[
\X_i^+:=\{x \in \X : S_N(x)=i \} \,\,\, \text{and} \,\,\,
\X_i^-:=\{x \in \X : S_N(x)=-i \}.
\]

The new proposal chain will be
\begin{equation}\label{proposal}
\begin{split}
K(x,y)&= p_1 K_E(x,y) +(1-p_1)K_0(x,y) \qquad
\text{if $x \in \X_0$} \\
K(x,y)&= p_1 K_E(x,y)+p_2\J_{\{-x\}}(y) + (1-p_1-p_2)K_i(x,y)  \\
& \qquad\qquad \qquad \qquad \qquad \qquad   \qquad \qquad \,\,    \text{if $x \in \X_i$, $i \not=0$}\\
\end{split}
\end{equation}
where $p_1,p_2$ belong to $(0,1)$, $p_1+p_2 <1$, and
\[
K_{i}(x,y)=
\J_{\X_i^+}\{x\} K_{i}^+(x,y)+\J_{\X_i^-}\{x\} K_{i}^-(x,y) \qquad
(i \not=0).
\]
We shall assume that
 $K^\pm_{i}$ ($K_{0}$, respectively) are irreducible, symmetric and aperiodic
chains on $\X_i^\pm$ ( $\X_0$, respectively).

As a leading example we shall take
\begin{equation}\label{esempio}
\begin{split}
&K_{0}(x,y)=\frac{1}{{N \choose N/2}} \qquad y \in \X_0
\\ &
K_{i}^{\pm}(x,y)=\frac{1}{{N \choose (N-i)/2}} \qquad y \in
\X_i^\pm, \\
\end{split}
\end{equation}
that is: a realization of a chain $K^{\pm}_i$ ($K_0$, respectively) is
simply a sequence of independent uniform random sampling from
$\X_i^{\pm}$ ($\X_0$, respectively).

\begin{remark}
Note that (\ref{esempio}) is the $(n,k)$-Bose-Einstein distribution
with $n=(N+i)/2$ and $k=(N-i)/2+1$ and recall that  there is a very easy way to
directly
generate Bose-Einstein configurations. One may place $n$ balls
sequentially
into $k$ boxes, each time choosing a box with probability proportional
to its current content plus one. Starting from the empty configuration
this results in a Bose-Einstein distribution for every stage.
\end{remark}

Now let $M$ be the Metropolis chain
 defined by the transition
kernel (\ref{MChain}) with $K$ as in (\ref{proposal}), i.e. for every $x$ in $\CX_i^{\pm}$ $(i\not=0)$
\[
M(x,y)=
\left\{
\begin{array}{ll}
\frac{p_1}{N}\min\left(1,\frac{\pi(y)}{\pi(x)}\right)&\mbox{ if } y=x^{(j)},\quad j=1...N\\
&\\
p_2& \mbox{ if } y=-x\\
&\\
(1-p_1-p_2) K_i^\pm(x,y) & \mbox{ if }  y \in \CX_i^{\pm},  y \not =x\\
&\\
1-\sum_{z\neq x}M(x,z)&\mbox{ if } y=x\\
\end{array}
\right.
\]
while for $x$ in $\CX_0$
\[
M(x,y)=
\left\{
\begin{array}{ll}
\frac{p_1}{N}\min\left(1,\frac{\pi(y)}{\pi(x)}\right)&\mbox{ if } y=x^{(j)},\quad j=1...N\\
&\\
(1-p_1)K_0(x,y)& \mbox{ if }  y \in \CX_0, y \not =x
\\
&\\
1-\sum_{z\neq x}M(x,z)&\mbox{ if } y=x.\\
\end{array}
\right.
\]
By construction
$M$ is an aperiodic, irreducible and reversible
chain with stationary distribution $\pi$.
Then, when (\ref{esempio}) holds true,
\[
M(x,y)=
\left\{
\begin{array}{ll}
\frac{p_1}{N}\min\left(1,\frac{\pi(y)}{\pi(x)}\right)&\mbox{ if } y=x^{(j)},\quad j=1...N\\
&\\
p_2& \mbox{ if } y=-x\\
&\\
(1-p_1-p_2)\frac{1}{{N \choose (N-i)/2}}& \mbox{ if }  y \in \CX_i^{\pm},  y \not =x\\
&\\
1-\sum_{z\neq x}M(x,z)&\mbox{ if } y=x\\
\end{array}
\right.
\]
for $x$ in $\CX_i^{\pm}$ $(i \not= 0)$, while if $x$ belongs to $\CX_0$
\[
M(x,y)=
\left\{
\begin{array}{ll}
\frac{p_1}{N}\min\left(1,\frac{\pi(y)}{\pi(x)}\right)&\mbox{ if } y=x^{(j)},\quad j=1...N\\
&\\
(1-p_1)\frac{1}{{N \choose N/2}}& \mbox{ if }  y \in \CX_0, y \not =x
\\
&\\
1-\sum_{z\neq x}M(x,z)&\mbox{ if } y=x.\\
\end{array}
\right.
\]

In order to bound the spectral gap of $M$ we shall use the
decomposition theorem described in Subsection \ref{chaindec}.
To this end, for every $i=0,2,\dots,N$ and every $j \not =i$ set
\[
\bar P(i,j):= \frac{1}{2 \pi(\X_i)} \sum_{x \in \X_i} \sum_{y \in \X_j}
M(x,y) \pi(x)
\]
and
\[
\bar P(i,i):= 1- \sum_{j \not =i }\bar P(i,j).
\]
As already noted,  $\bar P$ is a reversible chain on $\{0,2,\dots,N\}$
with stationary distribution
\[
\bar \pi(i):=\pi(\X_i).
\]
Moreover define for every $i=0,2,\dots,N$ a chain on $\X_i$ setting
\[
P_{\X_i}(x,y):=M(x,y)+ \J_{x}(y)\left(\sum_{z \in \X_i^c} M(x,z)\right)
\]
where both $x$ and $y$ belong to $\CX_i$.
In the same way, define chains on $\X_i^+$ and $\X_i^-$ for
$i=2,\dots,N$ setting
\[
P_{\X_i^\pm}(x,y):=P_{\X_i}(x,y) \quad (y \not = x, \,\, x,y \in X_i^\pm)
\]
and
\[
P_{\X_i^\pm}(x,x):=1-\sum_{y \in \X_i^\pm y\not=x }P_{\X_i}(x,y).
\]
These chains are reversible on $\X_i$ ($\X_i^\pm$, respectively)
and have as stationary distributions
\[
\pi_{\X_i}(x):=\frac{\pi(x)}{\pi(\X_i)}=\frac{1}{|\CX_i|}
\quad \text{and} \quad
\pi_{\X_i^\pm}(x):=\frac{\pi_{\X_i} (x)}{\pi_{\X_i}(\X_i^\pm)}
=\frac{1}{|\CX_i^{\pm}|},
\]
respectively.
Finally, for every $i=2,4,\dots,N$, define a chain on $\{+,- \}$
setting
\[
\begin{split}
P_i(+,-)&:= \frac{1}{2\pi_{\X_i}(\X_i^+) }  \sum_{x \in \X_i^+} \sum_{y
  \in \X_i^-} P_{\X_i}(x,y) \pi_{\X_i}(x) \\
P_i(-,+)&:= \frac{1}{2\pi_{\X_i}(\X_i^- )}  \sum_{x \in \X_i^-} \sum_{y
  \in \X_i^+} P_{\X_i}(x,y) \pi_{\X_i}(x). \\
\end{split}
\]
Now the lower bound (\ref{decLB}), applied two times yields
\begin{equation}\label{mainboundISING}
\begin{split}
Gap(M) & \geq \frac{1}{2} Gap(\bar P)
 \min_{i=0,2,\dots,N}\left
\{ Gap(P_{\X_i}) \right \} \\
&\geq \frac{1}{2} Gap(\bar P) \min \Big [ Gap(P_{\X_0}), \\
& \min_{i=2,\dots,N}\left
\{  \frac{1}{2}Gap(P_i) \min \{ Gap(P_{\X_i^+}),
Gap(P_{\X_i^-}) \} \right \}\Big ].\\
\end{split}
\end{equation}
Hence, to get a lower bound on $Gap(M)$  it is enough to obtain bounds
on the gaps of the chains $\bar P$, $P_{\X_0}$, $P_i$,
$P_{\X_i^\pm}$.

The most important of these bounds is given by the following

\begin{proposition}\label{gapbarP} $\bar P$ is a birth and death chain
on $\{0,2,\dots,N\}$, more precisely
\begin{equation}\label{chainbarP}
\begin{array}{ll}
\bar P(0,2)=\frac{p_1}{2} & \quad   \\
\bar P(i,i+2)=\frac{p_1}{4}\frac{N-i}{N} & \quad i \not=N,0  \\
\bar P(i,i-2)=\frac{p_1}{4}\frac{N+i}{N}
\exp\{ 2 \b (1-i)/N \} & \quad i \not=0 .
\end{array} 
\end{equation}
Moreover
\[
\lambda_1(\bar P) \leq 1-\frac{p_1}{16}\frac{1}{(N/2+1)^3}.
\]
and
\[
\lambda_{N/2}(\bar P) \geq  1- p_1.
\]
\end{proposition}

The proof of the previous proposition is based on a
bound for a birth and death chain, given in the Appendix, which can be of its own interest.

As for the others chains,
we have the following
\begin{lemma}\label{lemmucolo}
For every $i=2,4,\dots,N$
\[
\begin{split}
&Gap(P_{\X_i^\pm})\geq (1-p_1-p_2) Gap(K^\pm_{i})\\
&Gap(P_i)=p_2,\\
\end{split}
\]
moreover
\[
Gap(P_{\X_0})\geq (1-p_1) Gap(K_{0}).
\]
\end{lemma}

In this way, using (\ref{mainboundISING}),
  we can prove the main result of this section.

\begin{proposition}\label{prop1} Let $M$ be the Metropolis chain
  derived by the chain $K$ defined as in
(\ref{proposal}) then
\[
\begin{split}
Gap(M) & \geq \frac{p_1p_2}{32}\frac{1}{(N/2+1)^3}
\min  \Big [ \frac{(1-p_1)}{p_2} Gap( K_{0} ),\\ &   \frac{(1-p_1-p_2)}{2}
 \min_{i \not =0} \min\{ Gap(K_{i}^+),Gap(K_{i,}^-)\} \Big ].\\
\end{split}
\]
If $K_{i}^\pm$ and $K_0$ are defined as in (\ref{esempio})
then
\[
Gap(M) \geq \frac{p_1p_2}{32}   \frac{1}{(N/2+1)^3} \min\Big[
\frac{(1-p_1-p_2)}{2} , \frac{(1-p_1)}{p_2}\Big]
\]
for every $\b>0$ and $N \geq N_0$.
\end{proposition}

Proposition \ref{prop1} shows that the gap is polynomial in ${1}/{N}$
independently of $\b$. Hence, even when $\b >1$, the variance of the
metropolis estimate obtained with this proposal can not grow up faster
than a polynomial in $N$.

Note that if in Proposition \ref{prop1} we choose
\begin{equation}\label{pmod}
p_1=1-a/(2N), \quad p_2=a/N
\end{equation}
we get
\[
Gap(M) \geq  \frac{C}{N^5}.
\]
Hence, even with this choice, the Metropolis algorithm is still fast mixing for every $\beta$.
It is worth noticing that the mean computational cost of this Metropolis
 does not change with respect to the Metropolis
which uses the proposal $K_E$.
Indeed, in the case of the usual Metropolis, the computational cost
needed
 to go from $X_n$ to $X_{n+1}$
is $O(N)$, since it is essentially due to a sample of one number among
$N$ numbers
(we need to decide which coordinate to flip).
In the case of the "modified" proposal,
things are slight more complex. In this case, at the beginning,
 we have an extra
``toss''.
If with this fist toss we decide to flip at random a coordinate
the cost is still $O(N)$ but if we need to sample
from $K^\pm_{i}$ the cost is
$O(N^2)$ (in this last case we need to pick a sample from a Bose-Einstein
distribution). Hence, although
our algorithm is "sometime" more expensive, if we take $p_1$ and $p_2$ as
in (\ref{pmod}), we get that the mean cost of our algorithm is still $O(N)$.

\section{The mean--field Blume-Emery-Griffiths model}\label{S:BEG}
The Blume-Emery-Griffiths (BEG) model (see \cite{Blume}) is an important lattice--spin
model in statistical mechanics, it has been studied extensively as a
model of many diverse systems, including $He^3-He^4$ mixtures as
well solid--liquid--gas systems, microemulsions, semiconductor
alloys and electronic conduction models. See, for instance,
\cite{Blume,Sivardiere1,Sivardiere2,Sivardiere3,Newman83,Schick,Kivelson}.
We will focus our attention on a simplified mean--field version of
the BEG model. For a mathematical treatment of this mean--field
model see \cite{ellisBEG}. In what follows let  $\X:=\{-1,0,1 \}^N$,
$N$ being an even integer, and for every $\beta>0$ and $K>0$ let
$\pi_{\beta,K,N}$ be the probability defined by
\[
\pi(x)=\pi_{\beta,K,N}(x)=
\exp\{-\beta R_N(x)+\frac{K \beta}{N} S_N^2(x) \}Z_N^{-1}(\beta,K) \qquad (x \in \X)
\]
where
\[
Z_N(\beta,K)=Z_N:=\sum_{x\in\X} \exp \left \{-\beta R_N(x)+\frac{K \b}{N} S_N^2(x) \right \}
\]
is the normalization constant,
$$
S_N(x):=\sum_{i=1}^N x_i\qquad \text{and} \qquad
R_N(x):=\sum_{i=1}^N x_i^2\qquad x=(x_1,x_2,...,x_N).
$$

A natural Metropolis algorithm can be derived by
using the  proposal chain
\begin{equation}\label{localproposalBEG}
K_E(x,y)=\frac{1}{2N} \sum_{j=1}^N [\J_{\{x^{(+ j)} \}}(y)+\J_{\{x^{(- j)} \}}(y)]
\end{equation}
where $x^{(\pm j)}$ denotes the vector $(x_1,\dots, x_j \pm
1,\dots,x_N)$, with the convention that $2=-1$ and $-2=1$.

The next proposition shows that there exists a critical region of the
parameters space in which
the Metropolis chain is slowly mixing. More precisely, using
some results of \cite{ellisBEG}
it is quite straightforward to proove the following

\begin{proposition}\label{slowlyBEG} Ler $M_E$ be the Metropolis chain (with stationary
distribution $\pi$)
with
  proposal chain $K_E$ defined in {\rm(\ref{localproposalBEG})}.Then,
there exists a non decreasing function $\Gamma:(0,+\infty)
\to(0,+\infty)$
with $\lim_{x \to 0} \Gamma(x)=+\infty$
and $\lim_{x \to \infty} \Gamma(x)=\gamma_c \simeq  1.082$ such that
for every couple of positive parametrs $(\beta,K)$ with $K > \Gamma(\beta)$
\[
Gap(M_E) \leq C e^{-\Delta N}
\]
for suitable constants $C=C(\gamma,K)>0$ and $\Delta=\Delta(\gamma,K)>0$.
\end{proposition}

As in the case of the mean--field Ising model,
we intend to by pass the slowly mixing  problem of this Metropolis
chain by choosing a different proposal.
To understand which kind of proposal is reasonable, here we choose
\[
\CG=\CS_N \times \{+1,-1 \}
\]
with $\CG$ acting on  $\X=\{-1,0,1\}^N$ by
\[
g(x)=(e\cdot x_{\s(1)},\dots,e\cdot x_{\s(N)}) \qquad g=(\s,e).
\]
At this stage,
 decompose $\X$ as the union of its "energy sets", that is
$$
\X=\X_{0,0}\cup\X_{1,1}\cup\X_{0,2}\cup\X_{1,3}\cup X_{3,3}\cup ...\cup\X_{0,N}\cup\X_{2,N}\cup ...\X_{N,N}
$$
where
$$\X_{s,r}:=\{x\in X:|S_N|=s \mbox{ and } R_N(x)=r\}$$
$r=0, 1,2,...,N$ and $s=1,3,...,r$  if $r$ is odd and
$s=0,2,...,N$ if $r$ is even.
Moreover,
for $s=1,2,...,N$, set
$$
\X_{s,r}^+:=\{x\in X:S_N=s \mbox{ and } R_N(x)=r\}
$$
and
$$
\X_{s,r}^-:=\{x\in X: S_N=-s \mbox{ and }R_N(x)=r\}.
$$

Note again that
$O_x=\X_{s,r}$ with $s=S_N(x)$ and $r=R_N(x)$.
The new proposal chain will be
\begin{equation}\label{proposalbeg}
\begin{split}
K(x,y)&= p_1 K_E(x,y) +(1-p_1)K_{0,r}(x,y)  \qquad
\text{if $x \in \X_{0,r},\quad r=0,2,...,N$} \\
K(x,y)&= p_1 K_E(x,y)+p_2\J_{\{-x\}}(y) + (1-p_1-p_2)K_{s,r}(x,y)  \\ & \qquad \qquad\qquad
\qquad
\text{if $x \in \X_{s,r}, s \not =0  $} \\
\end{split}
\end{equation}
where $p_1,p_2$ belong to $(0,1)$, $p_1+p_2 <1$, and
\[
K_{s,r}(x,y)=
\J_{\X_{s,r}^+}\{x\} K_{s,r}^+(x,y)+\J_{\X_{s,r}^-}\{x\} K_{s,r}^-(x,y) \qquad
(s \not=0)
\]
%
with
\begin{equation}\label{esempio2}
\begin{split}
&K_{0,r}(x,y)=\frac{1}{{N \choose r}{r \choose r/2}} \qquad y \in \X_{0,r}
\\ &
K_{s,r}^{\pm}(x,y)=\frac{1}{{N \choose r}{r \choose (r-s)/2}} \qquad y \in
\X_{s,r}^\pm. \\
\end{split}
\end{equation}
Now let $M$ be the Metropolis chain
defined by the transition
kernel (\ref{MChain}) with $K$ as in (\ref{proposalbeg}), i.e. for
every
$x$ in $\CX_{s,r}^{\pm}$ $(s\not=0)$
\[
M(x,y)=
\left\{
\begin{array}{ll}
\frac{p_1}{2N}\min\left(1,\frac{\pi(y)}{\pi(x)}\right)&\mbox{ if } y=x^{(\pm j)},\quad j=1...N\\
&\\
p_2& \mbox{ if } y=-x\\
&\\
(1-p_1-p_2)\frac{1}{{N \choose r}{r \choose (r-s)/2}}& \mbox{ if }  y \in \CX_{s,r}^{\pm},  y \not =x\\
&\\
1-\sum_{z\neq x}M(x,z)&\mbox{ if } y=x,\\
\end{array}
\right.
\]
while if $x$ belongs to $\CX_{0,r}$
\[
M(x,y)=
\left\{
\begin{array}{ll}
\frac{p_1}{2N}\min\left(1,\frac{\pi(y)}{\pi(x)}\right)&\mbox{ if } y=x^{(\pm j)},\quad j=1...N\\
&\\
(1-p_1)\frac{1}{{N \choose r}{r \choose r/2}}& \mbox{ if }  y \in \CX_{0,r}, y \not =x
\\
&\\
1-\sum_{z\neq x}M(x,z)&\mbox{ if } y=x.\\
\end{array}
\right.
\]
By construction
$M$ is an aperiodic, irreducible and reversible
chain with stationary distribution $\pi$.

Also in this case, to bound the spectral gap of $M$,
we shall use the chain decomposition tools.
Let
$$
\mathbb{D}_N=\{(0,0),(1,1),(0,2),(2,2),(1,3),
(3,3),(0,4),(2,4),(4,4),...,(0,N),(2,N),...,(N,N)\}
$$
and, for every couple  $(s,r),(\tilde{s},\tilde{r})$ in $
\mathbb{D}_N$, with $(s,r)\neq (\tilde{s},\tilde{r})$, let
\[
\bar P((s,r),(\tilde{s},\tilde{r})):= \frac{1}{2 \pi(\X_{s,r})} \sum_{x \in \X_{s,r}} \sum_{y \in \X_{\tilde{s},\tilde{r}}}
M(x,y) \pi(x)
\]
and
\[
\bar P((s,r),(s,r)):= 1- \sum_{(\tilde{s},\tilde{r})\not =(s,r) }\bar P((s,r),(\tilde{s},\tilde{r})).
\]
Once again, note that $\bar P$ is a reversible chain on $\mathbb{D}_N$
with stationary distribution
\[
\bar \pi(s,r):=\pi(\X_{s,r}).
\]
Moreover, for every $(s,r)$ in
$\mathbb{D}_N$,  define a chain on $\X_{s,r}$ setting
\[
P_{\X_{s,r}}(x,y):=M(x,y)+ \J_{x}(y)\left(\sum_{z \in \X_{s,r}^c} M(x,z)\right)
\]
where both $x$ and $y$ belong to $\CX_{s,r}$.
In the same way, define chains on $\X_{s,r}^+$ and $\X_{s,r}^-$ for
$(s,r)$ in $\mathbb{D}_N$, $s\neq 0$, setting
\[
P_{\X_{s,r}^\pm}(x,y):=P_{\X_{s,r}}(x,y) \quad (y \not = x, \,\, x,y \in X_{s,r}^\pm)
\]
and
\[
P_{\X_{s,r}^\pm}(x,x):=1-\sum_{y \in \X_{s,r}^\pm y\not=x }P_{\X_{s,r}}(x,y).
\]
These chains are reversible on $\X_{s,r}$ ($\X_{s,r}^\pm$, respectively)
and have as stationary distributions
\[
\pi_{\X_{s,r}}(x):=\frac{\pi(x)}{\pi(\X_{s,r})}=\frac{1}{|\CX_{s,r}|}
\quad \text{and} \quad
\pi_{\X_{s,r}^\pm}(x):=\frac{\pi_{\X_{s,r}} (x)}{\pi_{\X_{s,r}}(\X_{s,r}^\pm)}
=\frac{1}{|\CX_{s,r}^{\pm}|},
\]
respectively.
Finally, for every $(s,r)$ in $\mathbb{D}_N$, $s\neq 0$,
define a chain on $\{+,- \}$
setting
\[
\begin{split}
P_{s,r}(+,-)&:= \frac{1}{2\pi_{\X_{s,r}}(\X_{s,r}^+) }  \sum_{x \in \X_{s,r}^+} \sum_{y
  \in \X_{s,r}^-} P_{\X_{s,r}}(x,y) \pi_{\X_{s,r}}(x) \\
P_{s,r}(-,+)&:= \frac{1}{2\pi_{\X_{s,r}}(\X_{s,r}^- )}  \sum_{x \in \X_{s,r}^-} \sum_{y
  \in \X_{s,r}^+} P_{\X_{s,r}}(x,y) \pi_{\X_{s,r}}(x). \\
\end{split}
\]
At this stage,  the lower bound (\ref{decLB}), applied two times, yields
\begin{equation}\label{gapsBEG}
\begin{split}
Gap(M) & \geq \frac{1}{2} Gap(\bar P)
 \min_{(s,r)\in\mathbb{D}_N}\left
\{ Gap(P_{\X_{s,r}}) \right \} \\
&\geq \frac{1}{2} Gap(\bar P) \min \Big [ \min_{r=0,2,...,N} \left\{ Gap(P_{\X_{0,r}})\right\}, \\
& \min_{(s,r)\in\mathbb{D}_N,s\neq 0}\left
\{  \frac{1}{2}Gap(P_{s,r}) \min \{ Gap(P_{\X_{s,r}^+}),
Gap(P_{\X_{s,r}^-}) \} \right \}\Big ].\\
\end{split}
\end{equation}

To derive from the last bound a more explicit bound we need
 some preliminary work.
The first result we need is exactly the analogous of
Lemma \ref{lemmucolo}.

\begin{lemma}\label{lemmaBEG-gapvari}
Fore every $r=1,\dots,N$
\[
Gap(P_{\X_{0,r}})\geq(1-p_1)Gap(K_{0,r})=(1-p_1),
\]
moreover, for every $(s,r)$ in $\mathbb{D}_N$ with $s\not=0$,
\[
Gap(P_{\X_{s,r}}^\pm)\geq(1-p_1-p_2)Gap(K^\pm_{s,r})=(1-p_1-p_2).
\]
Finally, for every  $(s,r)$ in $\mathbb{D}_N$,
\[
Gap(P_{s,r})=p_2.
\]
\end{lemma}

Hence, (\ref{gapsBEG}) can be rewritten as
\begin{equation}\label{boundMBEGbis}
Gap(M) \geq Gap(\bar P) \frac{p_2}{2} \min\{(1-p_1)/2,(1-p_1-p_2)/2\}.
\end{equation}
It remains to bound $Gap(\bar P)$. Unfortunately the
the analogous of Proposition \ref{gapbarP} is not so simple,
hence we shall require an additional hypothesis.
In what follows let
\[
\begin{tabular}{lclr}
&&&\\
$q_{|\![N]\!|}(r)$&:=&$\displaystyle{N \choose r}e^{-\beta r}\left [ {r \choose \frac{r}{2}}+2\sum_{i=0}^{\frac{r}{2}-1}{r \choose i}e^{\frac{k\beta}{N}(r-2i)^2}\right]$& if
$r$ is even\\
&&&\\
$q_{|\![N]\!|}(r):$&=&$\displaystyle{N \choose r}e^{-\beta r}\left [2\sum_{i=0}^{\frac{r-1}{2}}{r \choose i}e^{\frac{k\beta}{N}(r-2i)^2}\right]$&if $r$ is odd\\
&&&\\
\end{tabular}
\]
$r=0,1,\dots,N$ and set 
$$
\mathcal{A}=\{\beta>0,K>0 : \exists N_0\mbox{ such that }
\forall N\geq N_0,\quad q_{|\![N]\!|} \mbox{ is unimodal}\}.
$$

\begin{lemma}\label{lemmaBEG3}
For every $(\beta,K)$ in $\mathcal{A}$
$$
Gap(\bar P)\geq \frac{C p_1^2}{N^6}
$$
for a suitable constant $C=C(\beta,K)$.
\end{lemma}

Under the same assumptions of the previous Lemma
 we can  state the main results of this section.

\begin{proposition}\label{propBEG}
For every $(\beta,K)$ in $\mathcal{A}$
$$
Gap(M)\geq \frac{\tilde C p_1^2}{N^6}
$$
for a suitable constant $\tilde C=\tilde C(\beta,K)$.
\end{proposition}

We conjecture that $Gap(\bar P)$ is polynomial in $N$ for every $(\beta,K)$
such that $\beta\not=\Gamma(K)$ (where $\Gamma$ is the function of
Proposition \ref{slowlyBEG}), but we are not able to prove this
conjecture. In point of fact we
conjecture that
$\RE^+\times\RE^+ \setminus \{(\beta,K): \Gamma(K)=\beta \} \subset \mathcal{A}$. We plotted
$q_{|\![N]\!|}$ for different $N$, $\beta$ and $K$, and these plotts  seem, at least, to confirm that
$\RE^+\times\RE^+ \setminus \{(\beta,K): |\Gamma(K)-\beta| \leq \eps \} \subset \mathcal{A}$
for a suitable small $\epsilon$.
In Figure 1 we show the graph of $q_{|\![N]\!|}$ for few different $N$, $\b$
and $K$.

\begin{figure}
\includegraphics [scale=0.3]{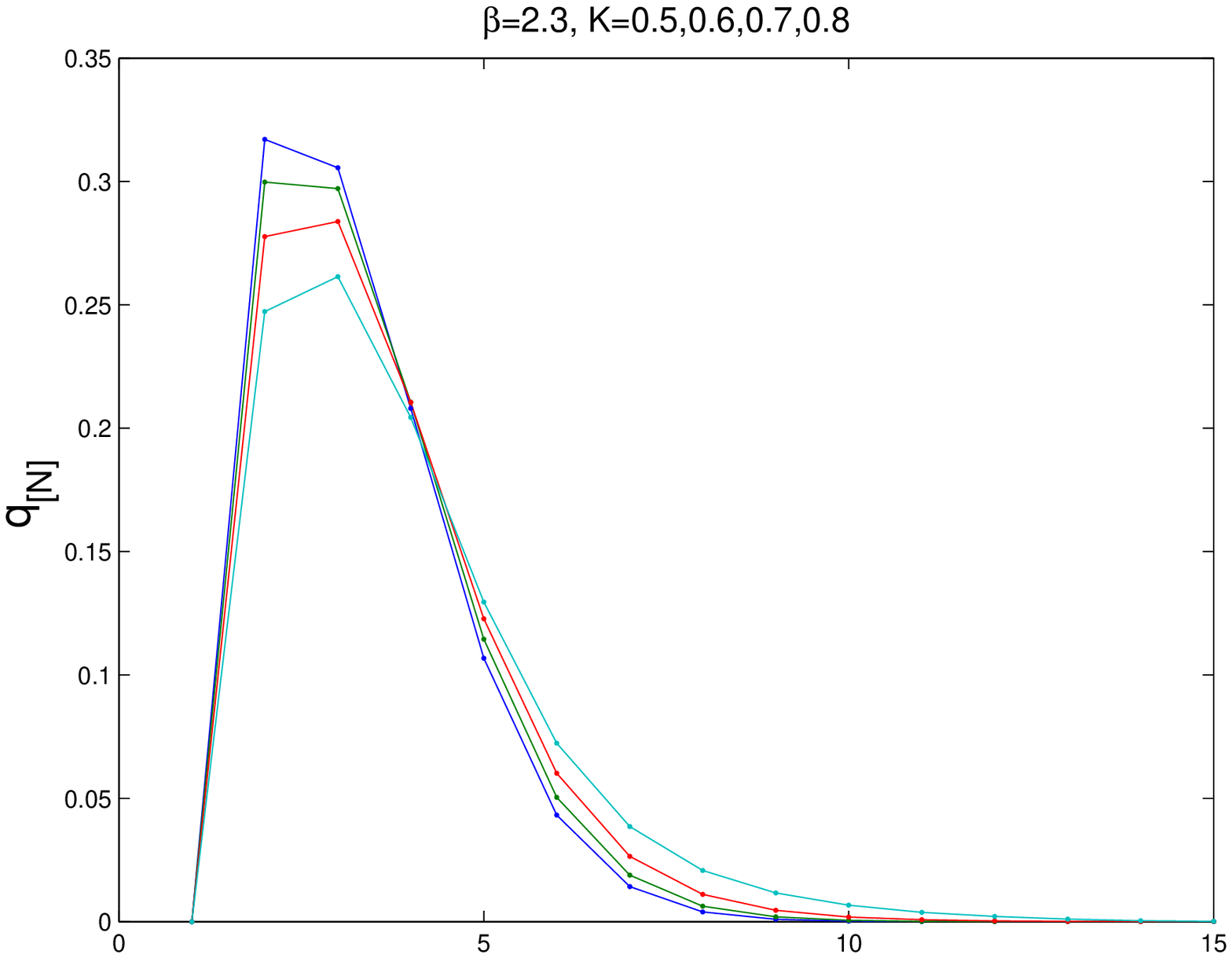}
\includegraphics [scale=0.3]{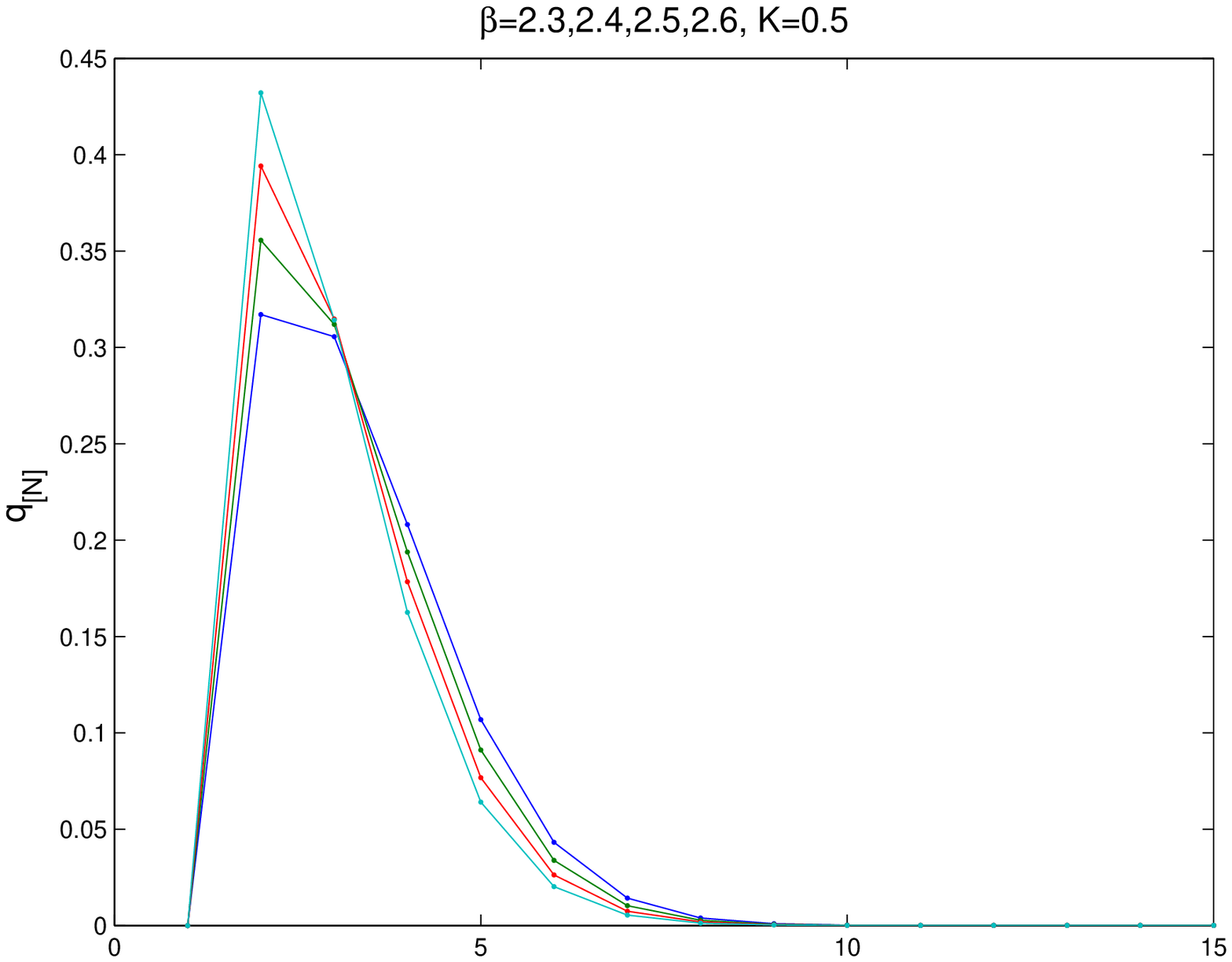}
\includegraphics [scale=0.3]{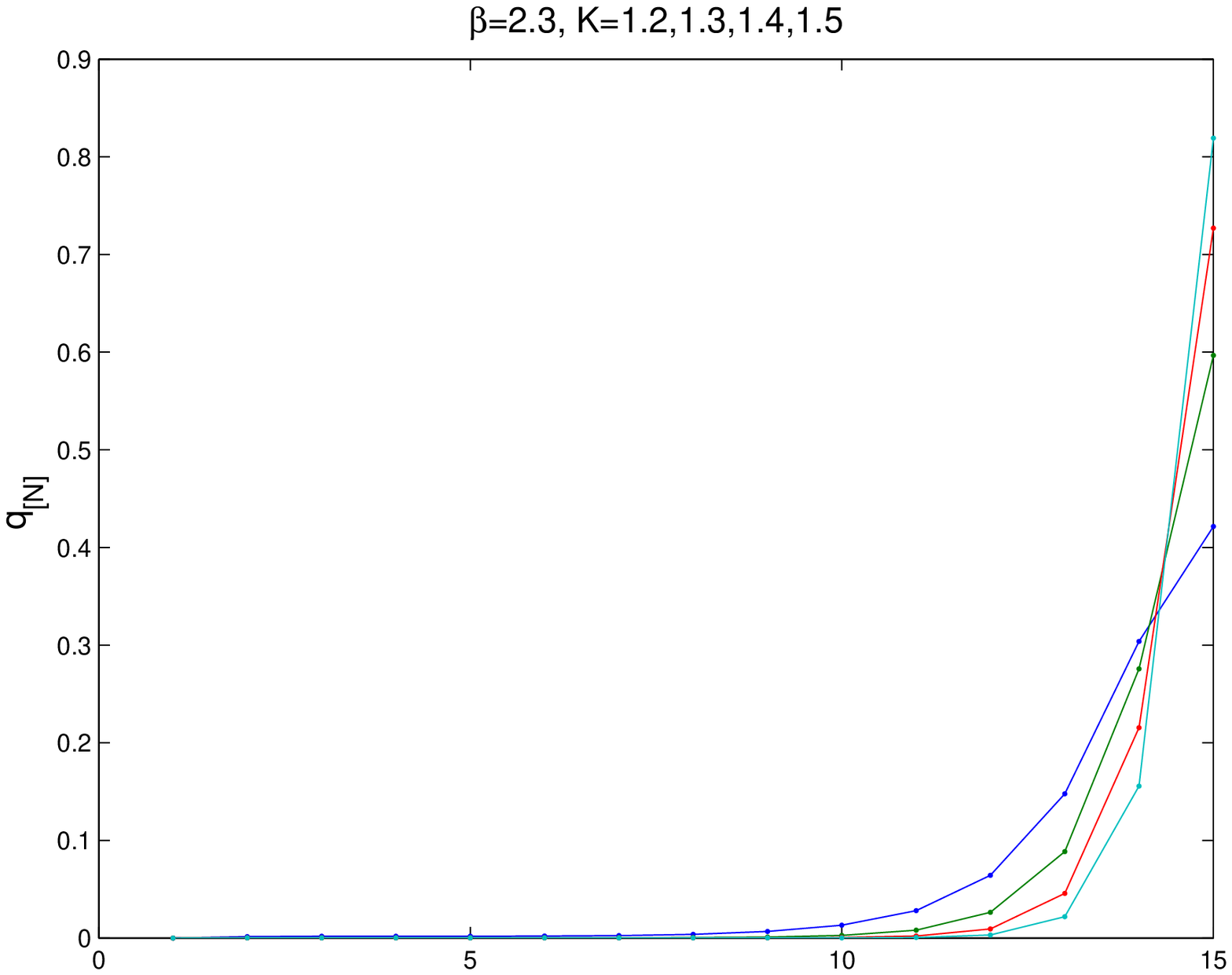}
\includegraphics [scale=0.3]{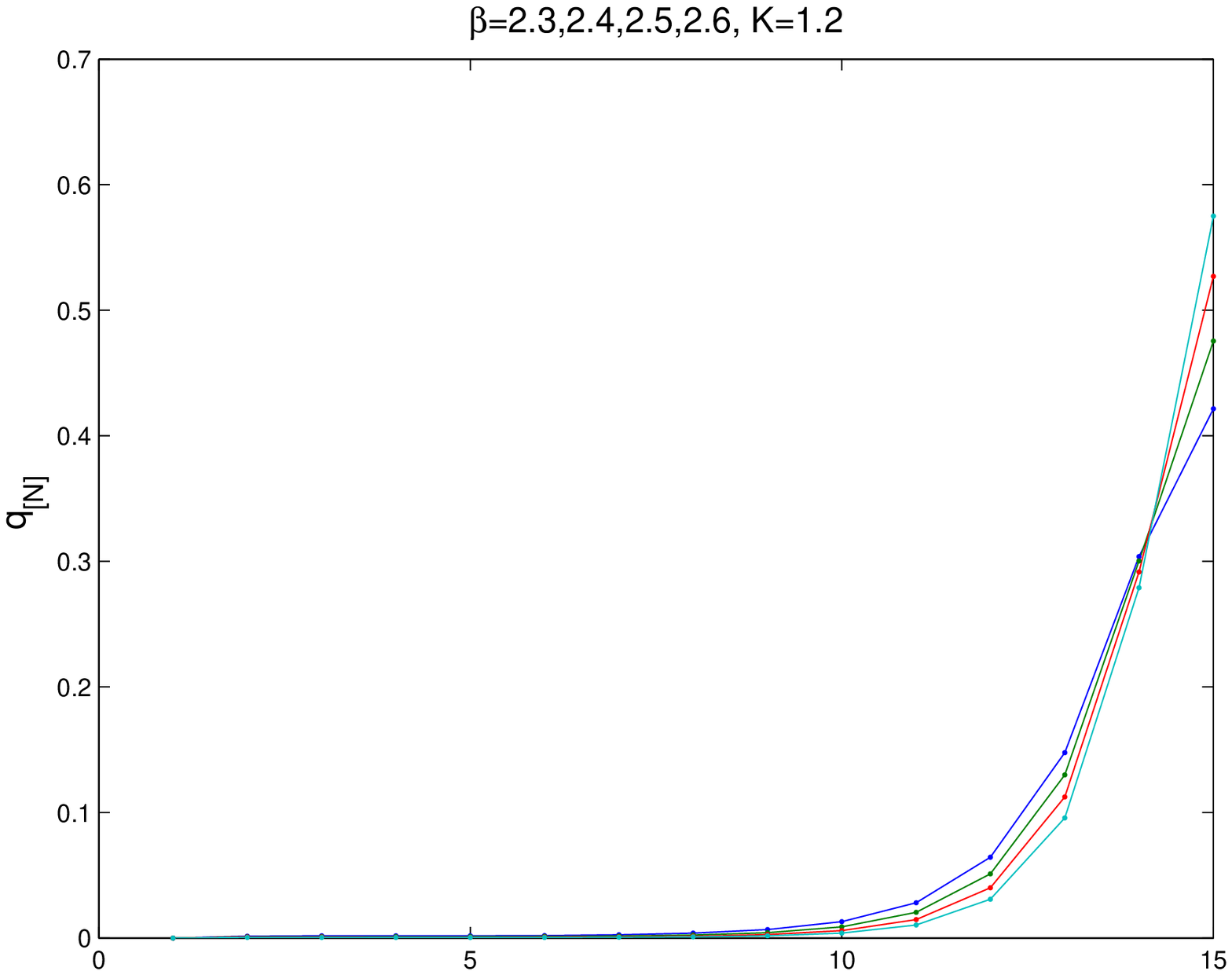}

     \caption{The function $q_{|\![N]\!|}$ for $N=15$ and few values of $\b$ and $K$.}
\end{figure}

\appendix
\section{The Spectral Gap of a Birth and Death Chain}
We derive here some bounds on the eigenvalues of
a birth and death chain that we shall use later.
These bounds are obtained using the so called geometric techniques, see
\cite{DiacStrook}.
Let $P_n$ be a birth and death chain on $\Omega_n=\{1,\dots,n
\}$. Assume that $P_n$ is reversible with respect to a probability
$p_n$, that is
\(
p_n(i)P_n(i,j)=p_n(j)P_n(j,i).
\)
Moreover let
\[
1> \lm_1 \geq \lm_2 \geq \dots \lm_{n-1} \geq -1
\]
the eigenvalues of $P_n$.

We can now prove the following variant of Proposition
6.3 in \cite{DSC2}.

\begin{lemma}\label{lemmaBD}
If there exist positive constants $A$, $q$, $B$ and  an integer $k$ such that
\[
\begin{split}
&P_n(i,i\pm1) \geq A n^{-q} \qquad (i \not = 1,n) \\
&P_n(1,2) \geq A n^{-q}  \\
&P_n(n,n-1) \geq A n^{-q}  \\
\end{split}
\]
and
\[
\begin{split}
p_n(i) & \leq B p_n(j) \quad i \leq j \leq k \\
p_n(j) & \leq B p_n(i) \quad  k \leq i \leq j \\
\end{split}
\]
then
\[
\lm_1 \leq 1- \frac{A}{B} \frac{1}{n^{q+2}}.
\]
\end{lemma}

\begin{flushleft}
\begin{proof}
We use the notation and the techniques of \cite{DiacStrook}, see also \cite{Bremaud} and \cite{DSC2}.
Choose the set of paths
$$
\Gamma=\{\gamma_{ij}=(i,i+1,...,j); i \leq j; i,j\in\Omega_n\}
$$
and for $e=(i,i+1)$ ($i < n$)
let
$$
\psi(e)=\frac{1}{p_n(i,i+1)}\sum_{\stackrel{\gamma_{l,m}\in\Gamma}{\gamma_{l,m}\ni
    e}}|\gamma_{l,m}|\frac{p_n(l)p_n(m)}{p_n(i)}
$$
where $|\gamma|$ is the length of the path $\gamma$. Setting
$K:=\sup_e \psi(e)$ one has
\[
\lambda_1 \leq 1- \frac{1}{K}
\]
(see Proposition 1' in \cite{DiacStrook}, or Exercise 6.4 page 248 in \cite{Bremaud}).
So, for our purposes,
 it suffices to give an upper bound on $K$.
Assume first that $e=(i,i+1)$
with $i < k \leq n$, since $|\gamma_{l,m}|\leq n$, it follows that
\begin{tabular}{lcl}
&&\\
$\psi(e)$&$\leq$&$\displaystyle\frac{n^q}{A}n\left(\sum_{\stackrel{s\geq i+1}{r\leq i}}\frac{p_n(r)p_n(s)}{p_n(i)}\right)$\\
&&\\
&$\leq$&$\displaystyle\frac{n^{q+1}}{A}\left(\sum_{r\leq i}\frac{p_n(r)}{p_n(i)}\right)\left(\sum_{s\geq i+1}p_n(s)\right)$\\
&&\\
&$\leq$&$\displaystyle\frac{n^{q+1}}{A}
\left(\sum_{r\leq i}B\right)\left(\sum_{s=1}^n p_n(s)\right)$\\
&&\\
&$\leq$&$\displaystyle n^{q+2}\frac{B}{A}.$\\
&&\\
\end{tabular}

All the other cases  can be treated in the same way.
Hence,
$$\sup_e \psi(e)\leq \frac{B}{A}n^{q+2}$$
and then
\[
\lm_1 \leq 1- \frac{A}{B} \frac{1}{n^{q+2}}.
\]
\end{proof}
\end{flushleft}
As for the smaller eigenvalues,
Gershgorin theorem yields  that
\[
\lm_{n-1} \geq -1 + 2 \min_i P(i,i).
\]
See, for instance, Corollary 2.1 in the Appendix of \cite{Bremaud}.
 Hence,
if there exists a positive constant $D$ such that
\[
P_n(i,i+1)+P_n(i,i-1) \leq D/2
\]
for every $i$, then
\begin{equation}\label{lowerboundBDchain}
\lm_{n-1} \geq 1 -D.
\end{equation}

\section{Proofs}

To prove Proposition \ref{gapbarP} we need first to show that $\bar \pi$ is essentially unimodal.

\begin{lemma}\label{lemmaunimodal}
Let
\[
q_N(i)= {N \choose \frac{N-i}{2}} \exp \left \{ \frac{\b}{2N} i^2  \right \}
\qquad i=0,2,4,\dots,N.
\]
For every $\b<1$ there exists an integer $N_0$ such that for every $N \geq N_0$
\[
q_N(i) \leq q_N(j)
\]
whenever $j \leq i$.
For every $\b \geq 1$ there exists an integer $N_0$ such that for every $N \geq N_0$
\[
q_N(i) \leq q_N(j)
\]
whenever $i \leq j \leq k_N$
and
\[
q_N(i) \geq q_N(j)
\]
whenever $k_N \leq i \leq j$, $k_N$ being a suitable integer.
\end{lemma}

\begin{proof}
Let $\Delta_N(i)$ be the ratio
$$
\Delta_N(i)=\frac{q_N(i+2)}{q_N(i)}\qquad i=0,2,4,...,N-2,
$$
so that
\[
\begin{split}
\Delta_N(i)&=\frac{{N \choose \frac{N-i}{2}-1}}{{N \choose
    \frac{N-i}{2}}}
\exp\left\{\frac{2\beta}{N}(1+i)\right\} \\ &=
\frac{N-i}{N+2+i}\exp\left\{\frac{2\beta}{N}(1+i)\right\}. \\
\end{split}
\]
Setting
$\Delta_N(x)=\frac{N-x}{N+2+x}\exp\left\{\frac{2\beta}{N}(1+x)\right\}$, $x$ in $[0,N-2]$,
it is enough to prove that $x \mapsto \Delta_N(x)$ takes the value $1$
at most once in  $[0,N-2]$, for sufficiently large $N$.
To prove this last claim first note that
\[
\begin{split}
\Delta_N(0)&=
\frac{N}{N+2}\exp\left\{\frac{2\beta}{N}\right\}
=\frac{1}{1+\frac{2}{N}}\exp\left\{\frac{2\beta}{N}\right\}\\
&=\left[1-\frac{2}{N}+2\left(\frac{2}{N}\right)^2+
o\left(\frac{1}{N^2}\right)\right]\left[1+\frac{2\beta}{N}
+\frac{1}{2}\left(\frac{2\beta}{N}\right)^2
+o\left(\frac{1}{N^2}\right)\right] \\
&=1-\frac{2}{N}(1-\beta)+\left(\frac{2}{N}\right)^2\left(\frac{\beta^2}{2}-\beta+2\right)+o\left(\frac{1}{N^2}
\right).\\
\end{split}
\]
Hence, there exists $N_0$ in $\mathbb{N}$ such that for $N\geq
N_0$:
\[
\begin{tabular}{lcl}
&&\\
$\beta\geq 1$&$\Rightarrow$&$\Delta_N(0)>1$\\
&&\\
$\beta<1$&$\Rightarrow$&$\Delta_N(0)<1$.\\
&&\\
\end{tabular}
\]
As for the first derivative note that
$$
\Delta_N'(x)=\frac{-2(N+1)+2\beta(N+2)-\frac{2\beta}{N}(x^2+2x)}{(N+x+2)^2}\exp\left\{\frac{2\beta}{N}(1+x)\right\},
$$
hence $\Delta_N'(x)=0$ if and only if
$$
-2(N+1)+2\beta(N+2)-\frac{2\beta}{N}(x^2+2x)=0.
$$
Rearranging the last equation as
$$
-\frac{2\beta}{N}x^2-\frac{4\beta}{N}+2[(\beta-1)N+2\beta-1]=0
$$
one sees that
the roots are
\[
x_{1,2}=1\pm\sqrt{1+\frac{2\beta-1}{\beta}N+\frac{\beta-1}{\beta}N^2}.
\]
Hence, after setting
\[
r:=1+\sqrt{1+\frac{2\beta-1}{\beta}N+\frac{\beta-1}{\beta}N^2} \qquad  \text{and} \qquad
\overline{r}:=1+\sqrt{1+N}
%
\]
one has
\[
\begin{tabular}{lcl}
&&\\
$\beta<1$&$\Rightarrow$&$\Delta_N'(x) < 0\qquad\forall x\in[0,N-2]$\\
&&\\
$\beta>1$&$\Rightarrow$&$\Delta_N'(x) > 0\qquad\mbox{for } x\in[0,r)$\\
&&$\Delta_N'(x)<0\qquad\mbox{for } x\in (r,N-2]$\\
&&\\
$\beta=1$&$\Rightarrow$&$\Delta_N'(x) <0\qquad\mbox{for } x\in[0,\overline{r})$\\
&&$\Delta_N'(x)<0\qquad\mbox{for } x\in (\overline{r},N-2]$\\
&&\\
\end{tabular}
\]
and this concludes the proof.
\end{proof}

\begin{proof}[Proof of Proposition \ref{gapbarP}]
By direct computations it is easy to prove (\ref{chainbarP}).
Hence
\[
P(i,i\pm2) \geq \frac{p_1}{4N}\geq\frac{p_1}{4(N+2)}\geq\frac{p_1}{8(\frac{N}{2}+1)},
\]
and
\[
P(i,i+2)+P(i,i-2) \leq \frac{p_1}{2}.
\]
Now observe that
\[
\bar \pi(0)=\frac{1}{Z_N(\b)}q_N(0)
\]
and
\[
\bar \pi(i)=\frac{2}{Z_N(\b)}q_N(i) \qquad i\not =0.
\]
Hence,
by Lemma \ref{lemmaunimodal}, if $\b<1$
\[
\bar \pi(i) \leq 2 \bar \pi(j)
\]
whenever $j \leq i$ and $N$ is large enough. While for
$\b > 1$
\[
\bar \pi(i) \leq  \bar \pi(j)
\]
whenever $i \leq j \leq k_N$
and
\[
\bar \pi(i) \geq \bar \pi(j)
\]
whenever $k_N \leq i \leq j$.
The thesis follows now by Lemma \ref{lemmaBD} and by (\ref{lowerboundBDchain}).
\end{proof}

In order to prove Lemma \ref{lemmucolo} we recall that
by Rayleigh's theorem
\begin{equation}\label{Rayleigh}
1-\lambda_1(P)= \inf \left \{\frac{\CE_p(f,f)}{Var_p(f)}: f
\,\,\text{nonconstant} \right\}
\end{equation}
where
\[
\CE_p(f,f):=<(I-P)f,f >_p=\frac{1}{2}\sum_{x,y}(f(x)-f(y))^2P(x,y)p(x),
\]
 $P$ being  a reversible chain
w.r.t. $p$, moreover
\begin{equation}\label{Rayleigh2}
1-|\lambda_{N-1}|= \inf \left \{\frac{\frac{1}{2}\sum_{x,y}(f(x)+f(y))^2P(x,y)p(x)}{Var_p(f)}: f
\,\,\text{nonconstant} \right\}
\end{equation}
(see, for instance, Theorem 2.3
 in Chapter 6 of \cite{Bremaud} and Section 2.1 of \cite{Dyeretal}).
%
At this stage set
\[
P_\eps(x,y):=(1-\eps)P(x,y)+\eps \J_{ x}(y).
\]
Hence, (\ref{Rayleigh}) yields
\[
\begin{split}
1-\lambda_1(P_\eps) &= \inf_{f \in L^2_p f \not = const}
\frac{\frac{1}{2} \sum_{x,y}(f(x)-f(y))^2P_\eps(x,y)p(x)  }{Var_p(f)} \\
 &= \inf_{f \in L^2_p f \not = const}
(1-\eps)
\frac{\frac{1}{2} \sum_{x \not =y}(f(x)-f(y))^2P(x,y)p(x)  }{Var_p(f)} \\
&=(1-\eps)(1-\lambda_1(P)).
\end{split}
\]
Arguing in the same way and using (\ref{Rayleigh2})
we get
\[
1-|\lambda_{|\CX|-1}(P_\eps)|\geq (1-\eps)(1-|\lambda_{|\CX|-1}(P)|).
\]
Hence,
\begin{equation}\label{gapeps}
Gap(P_\eps)\geq (1-\eps)Gap(P).
\end{equation}

\begin{proof}[Proof of Lemma \ref{lemmucolo}]
Note that
\[
P_{\X_i^\pm}(x,y)=(1-p_1-p_2)K^\pm_{i}(x,y) + (p_1+p_2)  \J_{x}(y)\\
\]
and, analogously,
\[
P_{\X_0}(x,y)=(1-p_1)K_{0}(x,y)+p_1 \J_{x}(y).\\
\]
Hence, by (\ref{gapeps}),
\[
Gap(P_{\X_i^\pm})\geq(1-p_1-p_2) Gap(K^\pm_{i})
\]
as well
\[
Gap(P_{\X_0})\geq(1-p_1) Gap(K_{0}).
\]
Finally note that $P_i$ is given by
\[
\left (
\begin{array}{cc}
1-\frac{p_2}{2} & \frac{p_2}{2} \\
\frac{p_2}{2} & 1-\frac{p_2}{2}
\end{array}
\right )
\]
for every $i$, hence $Gap(P_i)=p_2$.
\end{proof}

\begin{proof}[Proof of Proposition \ref{prop1}]  To prove the first part of the proposition it is enough to
combine Lemma \ref{lemmucolo},
 Proposition \ref{gapbarP} and (\ref{mainboundISING}).
To complete the proof observe that
$Gap(K_i^\pm)=Gap(K_0)=1$, when $K_i^{\pm}$ and $K_0$ are given by
(\ref{esempio}).
\end{proof}

In order to prove Proposition \ref{slowlyBEG} we
need some results obtained in \cite{ellisBEG}.

\begin{theorem}[Ellis-Otto-Touchette]\label{propositionLARGEDEVIATION}
Let $\rho_N$ be the distribution of $S_N(x)/N$ under $\pi_{\b,K,N}$,
then $\rho_N$ satisfies a large deviation principle on $[-1,1]$ with
rate function
\[
\tilde I_{\beta,K}(z)=J_\beta(z)  -\beta K z^2
-\inf_{t \in \RE}\{ J_\beta(t)-\beta Kt^2  \}
\]
with
\[
J_\beta(z)=\sup_{t \in \RE}\left\{tz - \log\left [ \frac{1+e^{-\b}(e^t+e^{-t}}{1+2e^{-\beta}} \right ] \right \}.
\]
Moreover, if $\tilde \CE_{\b,K}:=argmin \tilde I_{\beta,K}$,
then there exists a non decreasing function
$\Gamma :(0,+\infty) \to(0,+\infty)$
with $\lim_{x \to 0} \Gamma(x)=+\infty$
and $\lim_{x \to \infty} \Gamma(x)=\gamma_c \simeq  1.082$
such that for every $(\beta,K)$ with $K > \Gamma(\beta)$
then
\[
\tilde \CE_{\b,K}=\{\pm z(\beta,K) \not =0 \}.
\]
In particular, for such $(\b,K)$ and for every $0<\eps< |z(\beta,K)|$ there exists a constant
$C_1=C_1(\eps,\b,K)$ such that
\begin{equation}\label{ellis1}
\rho([0,\eps]) \leq  C_1 \exp\{ -\frac{N}{2} \gamma_{\eps,\beta,K} \}
\end{equation}
with
\begin{equation}\label{ellis2}
\gamma_{\eps,\beta,K}= \inf_{z \in [0,\eps] }  \tilde I_{\beta,K}(z)>0.
\end{equation}
\end{theorem}

\begin{proof} For the first part
see Theorems 3.3, 3.6 and 3.8 in \cite{ellisBEG}.
As for (\ref{ellis1})-(\ref{ellis2}), they
are standard consequences of the theory of the large deviations
and of the first part of the proposition, see, e.g., Proposition
6.4 of \cite{ellisLecture}.
\end{proof}

\begin{proof}[Proof of Proposition \ref{slowlyBEG}] We intend to use
  the Chegeer's inequality. To do this,
let $A:=\{x: S_N(x)<0 \}$, $B:=\{x: S_N(x)>0 \}$, $C:=\{S_N(x)=0 \}$. First of all
note that, by symmetry, $\pi(A)=\pi(B)=(1-\pi(C))/2 \leq 1/2$.
The main task is to bound
\[
\phi(A)= \sum_{x \in A} \sum_{y \in A^c} \pi(x) M_E(x,y)=\sum_{y \in A^c}
\sum_{x \in A} \pi(y) M_E(y,x).
\]
Now, observe that if $S_N(y)>1$ then $M_E(y,x)=0$ for every $x$ in $A$, hence
\[
\begin{split}
\phi(A)& = \sum_{y:S_N(y)=0 } \pi(y) \sum_{x \in A} M_E(y,x) +
\sum_{y:S_N(y)=1} \pi(y) \sum_{x \in A } M_E(y,x) \\
&\leq \pi \left \{y: S_N(y) \in \{0,1 \} \right \}. \\
\end{split}
\]
 This yields a bound on the conductance
\[
h= h(\pi,M_E) \leq \phi(A)/\pi(A) \leq \frac{ 2\pi \left\{ y: S_N(y) \in \{0,1 \} \right \} }{1- \pi\{y:
S_N(y)=0\}}.
\]
Now by Proposition \ref{propositionLARGEDEVIATION} we get
\[
h(\pi,M_E) \leq C_2 e^{-\Delta N}
\]
 for suitable constants $C_2$ and $\Delta>0$.
The thesis follows by Cheeger inequality  (\ref{cheger}).
\end{proof}

\begin{proof}[Proof of Lemma \ref{lemmaBEG-gapvari}]
The proof is exactly the same as the proof of Lemma \ref{lemmucolo}.
\end{proof}

In order to prove Lemma \ref{lemmaBEG3} it is convenient to fix some simple properties of the chain $\bar P$.

\begin{lemma}\label{lemmaBEGbarP} $\bar P$ is a random walk
on $\mathbb{D}_N$. If $\bar P((s,r),(\tilde s,\tilde r)) \not =0$,
\[
\bar P((s,r),(\tilde s,\tilde r)) \geq \frac{p_1 C_3}{N}
\]
for a suitable constant $C_3=C_3(\beta,K)$, moreover
\[
\bar P((s,r),(\tilde s,\tilde r)) \leq \frac{p_1}{4}
\]
for every $(s,r),(\tilde s, \tilde r))\not = ((0,0),(1,1))$.
\end{lemma}

\begin{proof}[Proof of Lemma \ref{lemmaBEGbarP}]
Easy but  tedious computations show that
\begin{equation*}\label{chainbarPbis}
\begin{split}
&\bar P((0,0),(1,1))=\frac{p_1}{2}\min\left(1,\exp\{\frac{K\beta}{N}-\beta\}\right)    \\
&\bar P((0,N),(1,N-1))=\frac{p_1}{4}   \\
&\bar P((0,N),(2,N))=\frac{p_1}{4}  \\
\end{split} 
\end{equation*}
\begin{equation*}
\begin{split}
&\bar P((0,r),(2,r))=\frac{p_1}{4N}\qquad r=0,2,4,...,N-2   \\
&\bar P((0,r),(1,r-1))=\frac{p_1}{4N} \qquad r=0,2,4,...,N-2 \\
&\bar P((0,r),(1,r+1))=\frac{p_1}{2N}\min\left(1,\exp\{\frac{K\beta}{N}-\beta\}\right)  \qquad r=0,2,4,...,N-2 \\
&\bar P((s,r),(s+2,r))=\frac{p_1}{8N}(r-s) \\
\end{split} 
\end{equation*}
\begin{equation*}
\begin{split}
&  \qquad \qquad (s,r)\in\mathbb{D}_N, 0<s\leq N-2, r\leq N  \\
&\bar P((s,r),(s-2,r))=\frac{p_1}{8N}(r+s)\exp\{4\frac{K\beta}{N}(1-s)\}  \\
& \qquad \qquad     (s,r)\in\mathbb{D}_N,0<s\leq N, r\leq N   \\
&\bar P((s,r),(s+1,r+1))=\frac{p_1}{4N}(N-r)\min\left(1,\exp\{\frac{K\beta}{N}(2s+1)-\beta\}\right)  \\
&   \qquad \qquad  (s,r)\in\mathbb{D}_N, 0<s,r\leq N-1, \\
&\bar P((s,r),(s-1,r+1))=\frac{p_1}{4N}(N-r)\exp\{\frac{K\beta}{N}(-2s+1)-\beta\}  \\
&   \qquad \qquad    (s,r)\in\mathbb{D}_N, 0<s,r\leq N-1, \\
&\bar P((s,r),(s+1,r-1))=\frac{p_1}{8N}(r-s) \\
\end{split} 
\end{equation*}
\begin{equation*}
\begin{split}
&\qquad \qquad   (s,r)\in\mathbb{D}_N, 0<r\leq N, 0<s\leq N-2 \\
& \bar P((s,r),(s-1,r-1))=\frac{p_1}{8N}(r+s)\min\left(1,\exp\{\frac{K\beta}{N}(2s+1)-\beta\}\right) \\
& \qquad \qquad     (s,r)\in\mathbb{D}_N, 0<r\leq N,0<s\leq r. \\
\end{split} 
\end{equation*}
At this stage the statement follows easily.
\end{proof}

\begin{proof}[Proof of Lemma \ref{lemmaBEG3}]
In order to obtain a bound on the gap of $\bar P$ we shall apply
another time the decomposition technique.
Write
$$
\mathbb{D}_N=\bar\X_1\cup\bar\X_2\cup\bar\X_3\cup...\cup\bar\X_N,
$$
where
$$
\bar\X_1=\{(0,0),(1,1)\}\qquad\bar\X_r=\{(u_1,u_2)\in\mathbb{D}_n: u_2=r\}.
$$
On $|\![N]\!|:=\{1,...,N\}$ define a chain $P_{|\![N]\!|}$ setting
\[
P_{|\![N]\!|}(i,j):= \frac{1}{2 \bar \pi(\bar \X_i)} \sum_{ a
\in \bar\X_i} \sum_{b \in \bar\X_j}
\bar P (a,b) \bar \pi(a)
\]
and
\[
P_{|\![N]\!|}(i,i):= 1- \sum_{j \not =i } P_{|\![N]\!|}(i,j).
\]
Again $P_{|\![N]\!|}$ is a reversible chain on $|\![N]\!|$
with stationary distribution
\[
\bar \pi_{|\![N]\!|}(i):=\bar \pi(\bar \X_i).
\]
Finally for every $r=1,2,\dots,N$ we define a chain on $\bar\X_r$ by setting
\[
P_{\bar\X_r}(a,b):=\bar P(a,b)+ \J_{a}(b)\left(\sum_{z \in \bar\X_r^c}
\bar P(a,z)\right)
\]
where both $a$ and $b$ belong to $\bar \X_r$.
Now note that for every $r=2,3,\dots,N$
$P_{\bar\X_r}$ is a birth and death chain on the state space
$\{(1,r),(3,r),\dots,(r,r) \}$ for $r$ odd and
$\{(0,r),(2,r),\dots,(r,r) \}$ for $r$ even.
Let
\[
q_r(s):= {r \choose (r-s)/2} e^{ \frac{\beta K}{N} s^2 }
\]
and, for $r$ even,
\[
q_r(0):= 2{r \choose r/2}.
\]
Now  observe that $P_{\bar\X_r}$
has stationary distribution
\[
\pi_{r}(s) \propto q_r(s)
\]
with $s=0,2,\dots,r$ if $r$ is even and
 $s=1,3,\dots,r$ if $r$ is odd.
First of all let $r\not=1$, by Lemma \ref{lemmaBEGbarP}  and Lemma
\ref{lemmaunimodal}, it is easy to check that
$(P_{\bar\X_r},\pi_{r})$
meets the condition of Lemma \ref{lemmaBD}
with
\[
B=2, \quad n=[(r+2)/2], \quad A=C_3 p_1[(r+2)/2] N^{-1}
\]
($[x]$ being the integer part of $x$) and then
\[
1-\lambda_1(P_{\bar\X_r}) \geq \frac{C_3 p_1[(r+2)/2]}{2N[(r+2)/2]^3}
\geq \frac{C_3 p_1}{2N^3}.
\]
Finally, Lemma \ref{lemmaBEGbarP} with
 (\ref{lowerboundBDchain}) yields
\[
\lambda_{|\bar\X_r|-1 }(P_{\bar\X_r}) \geq 1- p_1.
\]
Hence, for every $r \not=1$, we have proved that
\begin{equation}\label{gapsPBEGs}
Gap(P_{\bar\X_r}) \geq C_3/2 p_1 N^{-3}.
\end{equation}
For $r=1$
\[
P_{\bar\X_1}=
\left (
\begin{array}{ll}
1-\a_1/2 & \a_1/2 \\
\a_2/2 & 1-\a_2/2
\end{array}
\right )
\]
where
\[
\a_1:=\frac{p_1}{2N}
 \min \left(1,\exp\{\frac{3K\beta}{N}-\beta\}\right)
\]
\[
\a_2:=p_1 \min \left(1,\exp\{\frac{K\beta}{N}-\beta\}\right)
\]
So
\[
Gap(P_{\bar\X_1}) \geq 1-|\frac{2-\alpha_1-\alpha_2}{2}|=\frac{\alpha_1+\alpha_2}{2}
\]
where the last equality follows from the fact that $\frac{\alpha_1}{2}\leq \frac{1}{2}$ and $\frac{\alpha_2}{2}\leq \frac{1}{2}$.
Hence, for sufficiently large  $N$, it's easy to see that
\begin{equation}\label{gapsPBEG1}
Gap(P_{\bar\X_1}) \geq C_4 p_1 N^{-3}
\end{equation}
with $C_4=C_4(\beta,K)$.
At this stage (\ref{gapsPBEGs}) with (\ref{gapsPBEG1}) gives
\begin{equation}\label{gapsPBEGsfin}
Gap(P_{\bar\X_r}) \geq C_5 p_1 N^{-3}
\end{equation}
for all $r\in |\![N]\!|$.
As for the gap of $P_{|\![N]\!|}$, first of all note that $P_{|\![N]\!|}$ is a birth and death chain on $|\![N]\!|$.
From Lemma \ref{lemmaBEGbarP}
\[
P_{|\![N]\!|}(i,i+1):= \frac{1}{2 \bar \pi(\bar \X_i)} \sum_{ a
\in \bar\X_i} \sum_{b \in \bar\X_{i+1}}
\bar P (a,b) \bar \pi(a)\geq\frac{p_1 C_3}{N}\frac{1}{2 \bar \pi(\bar \X_i)}\sum_{ a
\in \bar\X_i} \sum_{b \in \bar\X_{i+1}}\bar \pi(a)\geq \frac{p_1C_3}{2N}
\]
and analogously,
\[
P_{|\![N]\!|}(i,i-1)\geq \frac{p_1C_3}{2N}.
\]

Now, for $r\not=1$
\[
\bar \pi_{|\![N]\!|}(r)=q_{|\![N]\!|}(r)/(\sum_{i=0}^N q_{|\![N]\!|}(i))
\]
while
\[
\bar \pi_{|\![N]\!|}(1)=(q_{|\![N]\!|}(1)+q_{|\![N]\!|}(0))/(\sum_{i=0}^N q_{|\![N]\!|}(i)).
\]

So, using the unimodality of $q_{|\![N]\!|}$, we can apply Lemma \ref{lemmaBEGbarP} with
$$
A=\frac{p_1C_3}{2} \qquad B=\frac{e^{-2\beta}}{2}
$$
which gives
$$
\lambda_1(P_{|\![N]\!|})\leq 1-\frac{p_1C_3}{e^{-2\beta}}\frac{1}{N^3}\leq 1-\frac{p_1C_3}{N^3}.
$$
Using another time Lemma
\ref{lemmaBEGbarP}, by (\ref{lowerboundBDchain}), we get
$$
\lambda_N(P_{|\![N]\!|})\geq 1-p_1.
$$
Combining this two bounds we  have
\begin{equation}\label{gapsPBEGSUN}
Gap(P_{|\![N]\!|}) \geq \frac{C_3p_1}{ N^3}
\end{equation}
and so from (\ref{decLB})
$$
Gap(\bar P)\geq \frac{Cp_1^2}{2N^6},
$$
$C$ being a suitable constant that depends by $\beta,K,C_3,C_4,C_5$.
\end{proof}

\begin{proof}[Proof of Proposition \ref{propBEG}]
Combine Lemma \ref{lemmaBEG3} with  \ref{boundMBEGbis}.
\end{proof}

\section*{Acknowledgments}
We should like to thank  Persi Diaconis
for useful discussions and for having
encouraged us during this work,
 Antonietta Mira for suggesting some
interesting references and Claudio Giberti
for helping to  improve
an earlier version of the paper.



\end{document}